\documentclass [10pt]{article}
\usepackage{amsmath}
\usepackage[left=3.9cm, right = 3.9cm]{geometry}

\usepackage{amsfonts, amssymb}

\usepackage[utf8]{inputenc} 

\def\be#1\ee{\begin{equation}#1\end{equation}}

\newenvironment{proof}[1][] {\noindent {\bf Proof#1:} }{\hspace*{\fill}$\square$\medskip\par}

\newtheorem{thm}{Theorem}
\newtheorem{lem}[thm]{Lemma}

\newtheorem{cor}[thm]{Corollary}

\def\P{{\mathbb{P}}}
\def\R{\mathbb{R}}
\def\E{\mathbb{E}\,}

\def\Z{{\mathbb Z}}

\def\AAA{{\mathcal{E}}}

\def\EE{{\mathcal E}}

\newcommand{\eps}{\varepsilon}

\def\HH{{\mathcal H}}

\def\tX{{\widetilde{X}}}
\def\tY{{\widetilde{Y}}}

\def \=L{\ {\buildrel\hbox{\scriptsize d }\over =}\ }

\begin{document}

\title{A limit theorem for the last exit time over a moving nonlinear boundary for a Gaussian process
}

\author{
	N.Karagodin\footnote{{\tt email: nikitus20@gmail.com}.}

}

\date{}

\maketitle

\begin{abstract} 
	
We prove a limit theorem on the convergence
of the distributions of the scaled last exit time over a slowly moving nonlinear boundary for a class of Gaussian stationary processes.
The limit is a double exponential (Gumbel) distribution.

\end{abstract}


\medskip
\bigskip

\bigskip

\section{Introduction}

Consider a stationary Gaussian process with continuous trajectories and its  
"last exit time over a moving boundary", i.e. the last instant when the process 
hits a boundary $\eps f(t)$, where $t$ denotes time and $\eps>0$ is a drift (or trend) parameter.
After this instant, the process stays forever {\it under} the boundary. We are interested 
in the asymptotic distribution of the last exit time when the parameter $\eps$ goes to zero. 
In this work, we prove a limit theorem on the convergence of the distribution of 
the properly centered and scaled last exit time to a double exponential (Gumbel) law.

A special case of this problem, for a particular process and a linear boundary, emerged in recent works 
\cite{ABL19,ABL20}, that provide a mathematical study of a physical model (Brownian chain break). In the work \cite{KarLif} the same question for a wide class of processes was studied. In this work we give a natural generalization of this result for a variety of boundaries.

As far as we know, the problem setting handling small trends is new, although the last 
exit time is a sufficiently popular object in the problems of economical applications; such 
as studies of ruin probabilities. In risk theory, for a centered Gaussian process with continuous trajectories $Y(t)$, the process
\[
R(t) = u + \eps f(t) - Y(t)
\]
represents company balance. For instance, when $f(t) = t$, this process can be used as the simplest model of an insurance company balance with a starting balance $u$, fixed income per time $\eps$ and stochastic expenses Y. In this setting the moments $\inf\{t: R(t)<0\}$ and $\max\{t : R(t)\leq 0\}$ are called the ruin time and the ultimate recovery time respectively. 

There are plenty of works dedicated to the moment of the process with trend hitting some level. In order to achieve positive result one has to balance between the variety of trends and that of the considered class of processes.

For instance, the classical ruin time is well studied. The result for a locally stationary Gaussian process $Y(t)$ is given in \cite{HusPit}.

In \cite{Deb} the asymptotic behaviour of the distance between the ruin time and the ultimate recovery time when $u$ goes to infinity and $Y(t)$ has stationary increments is studied.

In \cite{Hash} the Paris ruin time $\inf \{t: \forall \delta \in [0, \Delta]\; R(t-\delta)\leq0\}$ for a fixed $\Delta$ is considered.

On the other hand, in \cite{Sal} the ruin time and the recovery time for any smooth trend are studied for the standard Brownian motion $Y(t)$.

In those settings, however, as a rule, one considers processes with stationary increments and a fixed trend, see also \cite{Hus,ParRab}. 
In this work we consider the ultimate recovery time for a stationary process without starting balance and small ratio of trend to volatility. The covariance function of the process is of Hölder type at 0 and decreasing at infinity faster than $1/\ln t$. Moreover, the class of trends is quite wide, covering, among others, the cases $\eps t^\beta$ and $\eps \exp\{t^q\}$.

\section{Main result}

Throughout the work we use the following notation for positive functions $f(x)$ and $g(x)$ 
\begin{itemize}
	\item $f(x) \prec g(x)$ if $\exists C>0: f(x)< Cg(x)$ asymptotically as x approaches $x_0$,
	\item $f(x) \asymp g(x)$ if simultaneously $f(x)\succ g(x)$ and $f(x) \prec g(x)$,
	\item $f(x) \sim g(x)$ if $\lim_{x\to x_0} \tfrac fg = 1$.
\end{itemize}

Sometimes we consider $f_a(x)$ and $g_a(x)$ depending on a parameter $a$. We say that the first relation is uniform over $a \in K$ if the constant $C$ does not depend on $a$. The second relation is uniform over $a \in K$ if the convergence is uniform.

Let $Y(t),t\in \R$, be a real-valued centered stationary Gaussian process with covariance function $\rho(t):=\E [Y(t)Y(0)]$.
We make two assumptions on the covariance function: 
\\at zero (Pickands condition); for some $v>0, Q>0$, $\alpha\in(0,2]$,
\begin{eqnarray} \label{1assumption_0}
	\rho(t)=v^2(1-Q|t|^\alpha+o(|t|^\alpha)), \quad \textrm{as } t\to 0,
\end{eqnarray}
and at infinity (Berman condition);
\begin{eqnarray} \label{1assumption_infty}
	\rho(t)= o( (\ln t)^{-1} ),   \quad \textrm{as } t\to +\infty.
\end{eqnarray}

Berman condition appears in the context of limit
theorems for maxima of Gaussian stationary sequences and processes,  see \cite{Ber} and 
\cite[Theorem 3.8.2]{Gal}.

Define $Y$'s last exit time over a boundary $f$ as
\[
T= T(\eps):=\max\{t: Y(t)=\eps\,f(t)\}.
\]
We are interested in the asymptotic behaviour of $T(\eps)$ when $\eps$ goes to 0. That is why, if it is not specified, we consider all the asymptotic relations when $\eps$ goes to 0.

Clearly $T(\eps) = \max \{t : \frac{1}{v}Y(t) = \frac{\eps}{v} f(t)\} = T_1(\eps/v)$, where $T_1$ is the last exit time of the scaled process $\frac{1}{v}Y(t)$. Therefore, one can study only normalized processes with $v = 1$ and then substitute $\eps/v$ instead of $\eps$ in the result for generality. Let us denote $\eps_v = \eps / v$.
\medskip

There are 3 conditions on the boundary function $f$.
\begin{itemize}
	\item {\it Ultimate monotonicity}: when $x$ tends to infinity, $f(x)$ is strictly increasing, twice differentiable and tends to infinity.  
	\item {\it Restriction on growth rate}: for some $0<\lambda \le 1$ it is true that $f''(x)/f'(x)\asymp x^{-\lambda}$ when $x\to\infty$.
\end{itemize}

To describe the third restriction we need to introduce two parameters. Take $\gamma = \gamma(\eps)$ such that
\begin{equation}
	\label{gammadef}
	\gamma^2 = 2\ln\left[\frac{(f^{-1})'(1/\eps_v)}{\eps_v}\right] + o(1),\;\; \mbox{when }\eps\to0.
\end{equation}
and $\tau_0 = \tau_0(\eps)$ such that 
\[
f(\tau_0)=\frac{\gamma}{\eps_v}, \;\; \mbox{when }\eps\to 0.
\]
They are important, because, in fact, $\tau_0$ is the main term of $T(\eps)$ asymptotics and $\tau_0^{\lambda} \gamma^{-2}$ is the precision, where the stochastic part appears.
Assume that $f$ also satisfies the following condition
\begin{itemize}
	\item {\it Regularity}: for some $\kappa>0$, $\beta$ and $\tilde \beta$ the following properties hold when $\eps \to 0$
	\begin{equation}
		\label{reg1}
		(f^{-1})'(y/\eps_v)\sim y^{\beta} (f^{-1})'(1/\eps_v), \;\;\;\mbox{uniformly over } y\in [(1-\kappa)\gamma, (1+\kappa)\gamma],
	\end{equation}
	\begin{equation}
		\label{reg2}
		(f^{-1})'(y/\eps_v)= o\left(y^{\tilde \beta} (f^{-1})'(1/\eps_v)\right), \;\;\;\mbox{uniformly over } y\in [(1+\kappa)\gamma, \infty).
	\end{equation}
\end{itemize}
Then we have the following limit theorem. In this theorem the Pickands constant $\HH_{\alpha}$ appears, see \cite{Pick_const_review} for more information about it.

\begin{thm} \label{t:main} Let $Y(t),t\in \R$, be a real-valued centered stationary Gaussian process satisfying
	assumptions \eqref{1assumption_0} and \eqref{1assumption_infty}.  Let $f$ be a function, satisfying the conditions above.
	For the Pickands constant $\HH_\alpha$ define
	$c:=\frac{Q^{1/\alpha}\HH_\alpha}{\sqrt{2\pi}}$. Then for any $r\in \R$ it is true that
	\[
	\lim_{\eps_v \to 0} \P\left\{ \frac{T(\eps)-A_{\eps}}{B_{\eps}}\le r \right\} = \exp(-c\exp(-r)),
	\]
	for the shift and scaling constants defined via $\tau_0$ and $\gamma$ as follows
	\begin{eqnarray*}
		A_\eps&:=& \tau_0 + \frac{(2/\alpha + \beta - 2)\ln\gamma}{f'(\tau_0)\eps_v \gamma} + o\left(\frac{1}{f'(\tau_0)\eps_v\gamma}\right),
		\\
		B_\eps &:=& \frac{1+o(1)}{f'(\tau_0)\eps_v\gamma}.
	\end{eqnarray*}
\end{thm}
\medskip

Note that in the recent work \cite{KarLif} the linear boundary was considered for the same class of processes. The proof below is based on the same ideas. However, working with a much larger class of boundaries requires more technical analysis.

\medskip


\subsection{Examples}
	
	There are several interesting examples of boundaries, for which the normalizing constants can be found in an explicit form.
	
	\medskip
	
	\begin{cor}
		For a boundary function $f(x) = x^d, d>0$ and a proccess $Y$, constant $c$, parameter $\eps_v$ and the moment $T(\eps)$ defined as before one has
		\[
		\lim_{\eps\to 0}\P\left\{\frac{T(\eps) - A_\eps}{B_\eps} \leq r\right\} = \exp (-c \exp (-r)),
		\]
		where the shift and scaling constants can be found explicitly 
			\begin{eqnarray*}
			A_\eps &=& \frac{(-2\ln\eps_v)^{\frac{1}{2d}-1}}{d^{\frac{1}{2d}}\eps_v^{\frac{1}{d}}}
			\left(-2\ln\eps_v + \left(\frac{1}{\alpha}+\frac{1}{2d}-\frac{3}{2}\right)\ln(-2\ln\eps_v) \right.
			\\&-& \left. \left(\frac{1}{\alpha}+\frac{1}{2d}-\frac{1}{2}\right)\ln d \right) ,
			\\B_\eps &=& \frac{(-2\ln\eps_v)^{\frac{1}{2d}-1}(1+o(1))}{d^{\frac{1}{2d}}\eps_v^{\frac{1}{d}}}.
		\end{eqnarray*}	
	\end{cor}	
	\begin{proof}
	For $f(x)=x^d, d>0$, we get $f''(x)=(d-1) x^{-1} f'(x)$, i.e. $\lambda=1$. Moreover, $f^{-1}(x)=x^{1/d}$ and $(f^{-1})'(x) = \frac{1}{d}x^{1/d-1}$. Therefore, we take
	\[
	\gamma^2 = 2\ln \left(\frac{1}{d}\eps_v^{-1/d}\right)+o(1)=\frac{-2\ln\eps_v}{d}-2\ln d +o(1).
	\]
	In addition, for every $y$ it is true that
	\[
	(f^{-1})'(y/\eps_v)=\frac{1}{d}y^{1/d-1}\eps_v^{1-1/d}=y^{1/d-1}(f^{-1})'(1/\eps_v),
	\]
	i.e. the regularity condition with $\beta=\tilde\beta=1/d-1$ holds.
	\\Take
	\[
	\tau_0=f^{-1}(\gamma/\eps_v) = (\gamma/\eps_v)^{1/d}.
	\]
	Then
	\[
	f'(\tau_0) = d \frac{f(\tau_0)}{\tau_0} = d \frac{\gamma}{\eps_v\tau_0}.
	\]
	Hence, by substituting everything in Theorem $\ref{t:main}$ we obtain
	\[
	A_\eps = \tau_0 + \frac{(2/\alpha +\beta -2)\tau_0\ln\gamma}{\gamma^2} + o\left(\frac{\tau_0}{\gamma^2}\right); \;\;\;
	B_\eps = \frac{\tau_0+o(\tau_0)}{\gamma^2}.
	\]
	Transforming these expressions, we get
	\begin{eqnarray*}
		A_\eps &=& \frac{\gamma^{\frac{1}{d}}}{\eps_v^{\frac{1}{d}}}
		+\left(\frac{2}{\alpha}+\frac{1}{d}-3 \right)\frac{\gamma^{\frac{1}{d}-2}\ln\gamma}{\eps_v^{\frac{1}{d}}}+o\left(\frac{\gamma^{\frac{1}{d}-2}}{\eps_v^{\frac{1}{d}}}\right)
		\\ &=& \frac{(-2\ln\eps_v)^{\frac{1}{2d}-1}}{d^{\frac{1}{2d}}\eps_v^{\frac{1}{d}}}
		\left(-2\ln\eps_v+\left(\frac{1}{\alpha}+\frac{1}{2d}-\frac{3}{2}\right)
		\ln(-2\ln\eps_v) \right.
		\\&-& \left. \left(\frac{1}{\alpha}+\frac{1}{2d}-\frac{1}{2}\right)\ln d \right) ,
		\\B_\eps &=& \frac{\gamma^{\frac{1}{d}-2}}{\eps_v^{\frac{1}{d}}} = \frac{(-2\ln\eps_v)^{\frac{1}{2d}-1}(1+o(1))}{d^{\frac{1}{2d}}\eps_v^{\frac{1}{d}}}.
	\end{eqnarray*}	
	\end{proof}
	\bigskip
	
	\begin{cor}
		For a boundary function $f(x) = \exp\{x^q\}, 0<q<1$ and a proccess $Y(t)$, constant $c$, parameter $\eps_v$ and the moment $T(\eps)$ defined as before one has
		\[
		\lim_{\eps\to 0}\P\left\{\frac{T(\eps) - A_\eps}{B_\eps} \leq r\right\} = \exp (-c \exp (-r)),
		\]
		where the shift and scaling constants can be found explicitly 
		
		\begin{eqnarray*}
			A_\eps &=& \frac{(-\ln\eps_v)^{\frac{1}{q}-1}}{q} \left(q (-\ln\eps_v) + \frac{1}{2}\ln\ln(-\ln\eps_v)
			+ \frac{1}{2}\ln\left(\frac{2}{q}-2\right)\right.
			\\&+& \left.\left(\frac{1}{\alpha}-\frac{3}{2}\right)\frac{\ln\ln(-\ln\eps_v)}{(\frac{2}{q}-2)\ln(-\ln\eps_v)}
			+ \frac{\left(\frac{1}{\alpha}-\frac{3}{2}\right)\ln\left(\frac{2}{q}-2\right)-\ln q}{(\frac{2}{q}-2)\ln(-\ln\eps_v)}\right),
			\\B_\eps &=& \frac{(-\ln\eps_v)^{\frac{1}{q}-1}}{(2-2q)\ln(-\ln\eps_v)}.
		\end{eqnarray*}
	\end{cor}

	\begin{proof}	
	For $f(x)=\exp\{x^q\}, 0<q<1$ we have
	\[
	f'(x) = q x^{q-1} \exp\{x^q\}
	\]
	and
	\[
	f''(x) = q (q-1) x^{q-2} \exp\{x^q\}+ q^2 x^{2q-2} \exp\{x^q\}.
	\]
	Therefore, the restriction on growth rate holds, since
	\[
	f''(x)/f'(x) \sim q x^{-(1-q)}.
	\]
	In addition
	\[
	(f^{-1})'(t) = \frac{(\ln t)^{\frac{1}{q}-1}}{q t},
	\]
	Therefore,
	\[
	\gamma^2 = \left(\frac{2}{q}-2\right)\ln(-\ln\eps_v)-2\ln q + o(1).
	\]
	Then we get the regularity condition with $\beta = -1, \tilde \beta = 0$. Uniformly over $y\in [0.5 \gamma, 2\gamma]$, as $\eps\to 0$ we have
	\[
	(f^{-1})'(y/\eps_v) = \frac{(\ln y - \ln \eps_v)^{\frac{1}{q}-1}}{q y/\eps_v} \sim y^{-1} (f^{-1})'(1/\eps_v),
	\]
	and uniformly over $y \in [2\gamma, \infty]$, as $\eps \to 0$ we have
	\[
	(f^{-1})'(y/\eps_v) =  \frac{(\ln y - \ln \eps_v)^{\frac{1}{q}-1}}{q y/\eps_v} = o((f^{-1})'(1/\eps_v)).
	\]
	Now we can use Theorem \ref{t:main}. 
	Take
	\[
	\tau_0 = f^{-1}(\gamma/\eps_v) = (\ln \gamma - \ln \eps_v)^{\frac{1}{q}},
	\]
	then
	\[
	f'(\tau_0) \eps_v = q \tau_0^{q-1} f(\tau_0) \eps_v =q \tau_0^{q-1} \gamma.
	\]
	Therefore, in the formula for $A_\eps$ and $B_\eps$ the precision we are interested in is
	\[
	o\left(\frac{1}{f'(\tau_0)\eps_v\gamma}\right)=o\left(\frac{\tau_0^{1-q}}{q\gamma^2}\right) = o\left(\frac{(-\ln\eps_v)^{\frac{1}{q}-1}}{\gamma^2}\right).
	\]
	Consequently,
	\begin{eqnarray*}
		A_\eps = \tau_0 
		+ \frac{(2/\alpha+\beta-2)\tau_0^{1-q}\ln\gamma}{q \gamma^2} 
		+ o\left(\frac{\tau_0^{1-q}}{\gamma^2}\right);\;\;\;
		B_\eps = \frac{\tau_0^{1-q}(1+o(1))}{q\gamma^2}.
	\end{eqnarray*}	
	Let us write
	\[
	\tau_0 = (-\ln\eps_v)^{\frac{1}{q}}+\frac{1}{q}(-\ln\eps_v)^{\frac{1}{q}-1}\ln\gamma
	+ o\left(\frac{(-\ln\eps_v)^{\frac{1}{q}-1}}{\gamma^2}\right),
	\]
	\begin{eqnarray*}
		\ln\gamma &=& \frac{1}{2}\ln\gamma^2 = \frac{1}{2}\left(\ln\left(\left(\frac{2}{q}-2\right)\ln(-\ln\eps_v)\right)
		+\ln\left(1 - \frac{q}{(\frac{1}{q}-1)\ln(-\ln\eps_v)}\right)\right)
		\\ &=& \frac{1}{2}\ln\ln(-\ln\eps_v) + \frac{1}{2}\ln\left(\frac{2}{q}-2\right) - \frac{\ln q}{(\frac{2}{q}-2)\ln(-\ln\eps_v)} + o\left(\frac{1}{\ln(-\ln\eps_v)}\right).
	\end{eqnarray*}
	After some tranformations we get
	\begin{eqnarray*}
		A_\eps &=& \frac{(-\ln\eps_v)^{\frac{1}{q}-1}}{q} \left(q (-\ln\eps_v) + \frac{1}{2}\ln\ln(-\ln\eps_v)
		+ \frac{1}{2}\ln\left(\frac{2}{q}-2\right)\right.
		\\&+& \left.\left(\frac{1}{\alpha}-\frac{3}{2}\right)\frac{\ln\ln(-\ln\eps_v)}{(\frac{2}{q}-2)\ln(-\ln\eps_v)}
		+ \frac{\left(\frac{1}{\alpha}-\frac{3}{2}\right)\ln\left(\frac{2}{q}-2\right)-\ln q}{(\frac{2}{q}-2)\ln(-\ln\eps_v)}\right),
		\\B_\eps &=& \frac{(-\ln\eps_v)^{\frac{1}{q}-1}}{(2-2q)\ln(-\ln\eps_v)}.
	\end{eqnarray*}
	\end{proof}
\section{Proof of  the main result}

The first restriction on the covariance function \eqref{1assumption_0} appears in the following lemma that serves below as one out of two basic tools
in our calculations.

\begin{lem} \label{1lem:pikpit} (Pickands--Piterbarg lemma).
	Let $Y(t),t\in \R$, be a real-valued centered stationary Gaussian process satisfying Pickands condition
	\eqref{1assumption_0} and such that
	\[
	\limsup_{t\to\infty} \rho(t) <1.
	\]
	Then for a Pickands constant $\HH_\alpha$ (in particular, $\HH_1=1$,$\HH_2=\pi^{-1/2}$), the following is true.
	\[
	\P \left\{ \max_{s\in[0,t]} Y(s) \ge x \right\}
	\sim \frac{Q^{1/\alpha}\HH_\alpha}{\sqrt{2\pi}} \cdot t \cdot (x/v)^{2/\alpha-1} e^{-x^2/2v^2}
	\]
	for all $x$ and $t$ such that the right hand side tends to zero and $tx^{2/\alpha}\to\infty$. 
\end{lem}

A first version of this lemma with fixed $t$ was obtained by Pickands \cite{Pick}, while
this version with variable $t$ (which is very important for our results) is due to Piterbarg,
see \cite[Lecture 9, Theorem 9.3.1]{Pit}.
\medskip

\bigskip

As was mentioned before, by scaling $Y(t)= v{\widetilde Y}(t)$ one may reduce the problem to the case $v=1$. That is why from now on $\eps_v = \eps$.

Let us fix $r\in \R$ and let
\[
\tau=\tau(\eps,r):= A_\eps +B_\eps r.
\]

\bigskip

The theorem's statement is equivalent to
\[
\lim_{\eps\to 0} \P\left\{ T(\eps)> \tau \right\} = 1 - \exp(-c\exp(-r)).
\]
Basic plan of the proof is the following. We are interested in the event of crossing the boundary $\eps f(t)$ after the moment $\tau$. 

First of all, for appropriatly large $\sigma$ this crossing happens before the moment $\sigma$ with probability close to 1.

Now, to analyze the event of crossing somewhere in $[\tau, \sigma]$ we divide $[\tau, \sigma]$ into intervals and approximate the boundary $\eps f(t)$ with a staircase, i.e. a function that is constant on each interval and has jumps in between. The event of crossing the staircase is much simpler and we can use Lemma \ref{1lem:pikpit} to analyze it, by studying the events of crossing on different intervals. If this events of crossing were almost independent for different intervals, we would have to work with a sum of independent random variables. The only problem here is that the events of crossing on the neighboring intervals are unpredictably correlated. 

To fix it one has to choose intervals appropriatly. We divide $[\tau, \sigma]$ into alternating long and short intervals of length $\ell(\eps)$ and $s(\eps)$ correspondingly. Then we show that the crossing happens on one of the small intervals with probability close to zero. The role of the small intervals is to separate long intervals, so we can controll the correlation of the process between different long intervals via Berman condition \eqref{1assumption_infty}.

Finally, to work with the crossing on different long intervals we concentrate the correlation between different intervals in one auxilary term using Slepian inequality \eqref{1slepian}. This part of the proof starts at Lemma \ref{Slepian_original}. It allows us to pass to the sum of independent random variables and then find the probability we are interested in.

The interval lengths $\ell=\ell(\eps)$, $s=s(\eps)$ must satisfy the following relations
\begin{eqnarray} \label{ls0}
	\ln s \succ \gamma^2
	\\  \label{ls1}
	s/\ell\to 0,
	\\  \label{ls2}
	f(\tau) f'(\tau) \eps^2 \, \ell \to 0.
\end{eqnarray}
We need the first relation to get a lower bound on the correlation between different long intervals. The second relation implies that the probability of crossing on some small interval is  close to 0. The third one is equivalent to the fact that $\ell = o(\sigma - \tau)$, i.e. it means that the number of intervals grows to infinity and the staircase is a good enough approximation of the boundary.

One can find the proof that it is possible to pick $\ell$ and $s$ satisfying these conditions in Section \ref{ss:proof_prop}.

We cover the halfline $[\tau,\infty)$ with the following system of sets:
\begin{itemize}
	\item
	the halfline $[\sigma,\infty)$, where $\sigma:=A_\eps +B_\eps R$ and $R=R(\eps)$ slowly 
	tends to infinity, as $\eps\to 0$.
	The choice of $R$ is further specified at the end of the proof but note that we assume only upper bounds on the growth rate of $R$.
	\item
	long intervals $L_i=[(\ell+s)i,(\ell+s)i+\ell]$, $i\in\Z$, of length  $\ell = \ell(\eps)$,
	\item
	short intervals $S_i=[(\ell+s)i+\ell,(\ell+s)(i+1)]$, $i\in\Z$, of length $s = s(\eps)$.
	
\end{itemize}
	
Let
\[
X_i^\eps:= \max_{t\in L_i} Y(t); \quad  V_i^\eps:= \max_{t\in S_i} Y(t).
\]
By using stationarity, we infer from the Pickands--Piterbarg lemma the asymptotics
\[
\P\{ X_i^\eps\ge x\} \sim c\, \ell\, x^{2/\alpha-1} \exp(-x^2/2)
\]
and 
\[
\P\{ V_i^\eps\ge x\} \sim c\, s\, x^{2/\alpha-1} \exp(-x^2/2),
\]
as soon as the corresponding right hand sides tend to zero and $\ell x^{2/\alpha} \to \infty$.

Define the index sets
\begin{eqnarray*}
	I_1 &:=& \{i:\ \substack{(\ell+s)i+\ell\ge\tau \\(\ell+s)i < \sigma }\},
	\\
	I_2 &:=& \{i:\ \substack{(\ell+s)i \ge\tau \\(\ell+s)i+\ell < \sigma}\},
	\\
	I_3 &:=& \{ i: \substack{(\ell+s)(i+1) \ge\tau \\(\ell+s)i+\ell < \sigma}\}
\end{eqnarray*}
chosen so that the following inclusions hold:
\be \label{inclusions}
\bigcup_{i\in I_2}  L_i  \subset  [\tau,\sigma] \subset  \big( \bigcup_{i\in I_1} L_i \big) \cup \big( \bigcup_{i\in I_3} S_i \big).
\ee

Let us define the events related to the exits of our process over the boundary:
\begin{eqnarray*}
	\AAA_1 &:=& \bigcup_{i\in I_1} \{ X_i^\eps \ge f((\ell+s)i)\eps \},  
	\\
	\AAA_2 &:=&  \bigcup_{i\in I_2}  \{ X_i^\eps\ge f((\ell+s)(i+1))\eps \},  
	\\
	\AAA_3 &:=& \bigcup_{i\in I_3}   \{ V_i^\eps\ge f((\ell+s)i)\eps \},  
	\\
	\AAA_4 &:=&  \{ \exists\, t>\sigma : Y(t) \ge \eps f(t) \}.
\end{eqnarray*}

By using inclusions \eqref{inclusions} and ultimate monotonicity of $f$,
it is easy to see that the following inequalities are true:
\begin{eqnarray*}
	\P\{ T(\eps)>\tau \} &=& \P\{ \exists\, t> \tau: Y(t) \ge \eps f(t) \} \le \P\{\AAA_1\}+\P\{\AAA_3\}+\P\{\AAA_4\},
	\\
	\P\{ T(\eps)>\tau \} &\ge& \P\{\AAA_2\}.
\end{eqnarray*}
Therefore, it is sufficient to prove that, as $\eps\to 0$,  we have
\begin{eqnarray*}
	\P\{\AAA_1\}, \P\{\AAA_2\}  &\rightarrow& 1-\exp(-c\exp(-r)),
	\\
	\P\{\AAA_3\}, \P\{\AAA_4\} &\rightarrow& 0.
\end{eqnarray*}

The parameter $R$ should grow to infinity so slowly that uniformly over $\omega \in [-1, 1]$ we have
\begin{eqnarray} \label{sigma0}
	\tau_0 \sim \tau_0 + \omega(\sigma-\tau_0),\;\;
	f'(\tau_0) \sim f'(\tau_0 + \omega(\sigma-\tau_0)),\;\;
	f(\tau_0) \sim f(\tau_0 + \omega(\sigma-\tau_0)).
\end{eqnarray}
In other words, $f$ and $f'$ does not change a lot in the interval around $\tau_0$ containing $\sigma$. 

The following lemma provides some useful asymptotic estimates and shows that we can pick $R$ satisfying \eqref{sigma0}. Since it is only technical, the proof is postponed to Section \ref{ss:proof_prop}.
\begin{lem}
	\label{l:R}
	Consider any function $\tilde{R}(\eps)$ growing so slowly, that
	$\tilde{R} = o(\ln\gamma)$ and
	\[
	B_\eps \tilde{R} = \frac{\tilde{R}}{f'(\tau_0)\eps\gamma} + o\left(\frac{1}{f'(\tau_0)\eps\gamma}\right)
	\]
	Then for $A_\eps$ and $B_\eps$ described in Theorem $\ref{t:main}$ we have
	\begin{equation}
		\label{ABeq}
		f(A_\eps+ B_\eps \tilde{R})\eps = \gamma + \left(\frac{2}{\alpha}+\beta-2\right)\frac{\ln\gamma}{\gamma} + \frac{\tilde{R}}{\gamma} + o\left(\frac{1}{\gamma}\right).
	\end{equation}
	Besides that, for $R(\eps)$ corresponding to the same restriction, uniformly over $\omega \in [-1, 1]$ one has
	\[
	f'(\tau_0) \sim f'(\tau_0 + \omega (\sigma - \tau_0)), \;\;\; f(\tau_0) \sim f(\tau_0 + \omega (\sigma - \tau_0)).
	\]
	
\end{lem}

\bigskip

From the relation $\eqref{ABeq}$ applied to $\tilde R = r$ and $\tilde R = R$, we obtain
\begin{eqnarray} \label{tauassymp_1}
	f(\tau)\eps &\sim& \gamma ,
	\\ \label{tauassymp_0}
	(f(\tau)\eps)^{2/\alpha+\beta-2}\exp\{-(f(\tau)\eps)^2/2\}\}&\sim& e^{-r-\gamma^2/2},
	\\ \label{sigmaassymp_0}
	(f(\sigma)\eps)^{2/\alpha+\beta-2}\exp\{-(f(\sigma)\eps)^2/2\} &=& o\left(e^{-\gamma^2/2}\right).
\end{eqnarray}
These relations are crucial for our choice of $\tau$ and $\gamma$, because they appear in the resulting probabilities of $\EE_1, \EE_2, \EE_3, \EE_4$.

For the following relations we use an analogue of the restriction on growth rate condition but for $f$ instead of $f'$, namely the fact that
\begin{equation}
	\label{restr_growth}
\frac{f'(x)}{f(x)} \asymp x^{-\lambda},\;\;\; x\to \infty.
\end{equation}
To check this note that $(f'(x) x^{\lambda})' = f''(x) x^{\lambda} + \lambda f'(x) x^{\lambda - 1}$. Therefore, since $c f'(x) < f''(x)x^{\lambda} < C f'(x)$ and $\lambda \leq 1$, we can write a similar inequality but with new constants depending on $\lambda$
\[
c_\lambda f'(x) < (f'(x) x^{\lambda})' < C_\lambda f'(x), \;\;\; \mbox{for large enough x}.
\]
Now we integrate these inequalities and get
\[
c_\lambda f(x) + \textrm{const} < f'(x) x^{\lambda} < C_\lambda f(x) + \textrm{Const}.
\]
Since $f(x)\to \infty$ as $x\to\infty$, we divide by $f(x)$ and conclude that $f'(x) x^{\lambda} \asymp f(x)$, i.e. $f'(x)/f(x) \asymp x^{-\lambda}$.

Let us notice that due to the fact we proved \eqref{restr_growth} and obtained relation \eqref{tauassymp_1} we have
\begin{equation}
	\label{ff'_asymp} 
	f(\tau)f'(\tau)\eps^2\asymp (f(\tau)\eps)^2 \tau^{-\lambda} \asymp \gamma^2 \tau^{-\lambda}.
\end{equation}

Moreover, there is a connection between $\tau$ and $\gamma$.
Due to the relation \eqref{restr_growth} and regularity condition combined with definition of $\gamma$ \eqref{gammadef} we get
\begin{equation*}
(f^{-1}(\gamma/\eps))^{\lambda} \asymp \frac{f(f^{-1}(\gamma/\eps))}{f'(f^{-1}(\gamma/\eps))} = \frac{\gamma}{\eps} (f^{-1})'(\gamma/\eps) \sim \frac{\gamma^{1+\beta}}{\eps}(f^{-1})'(1/\eps) \sim \gamma^{1+\beta} e^{\gamma^2/2},
\end{equation*}
and since $\tau\sim\tau_0 = f^{-1}(\gamma/\eps)$, we get
\begin{equation}
	\label{tau_order}
	\tau^{\lambda} \asymp \gamma^{1+\beta} e^{\gamma^2/2}.
\end{equation}	

\bigskip

Now, after obtaining crucial relations, we proceed with the proof.
Let us first show that the probabilities of the events $\AAA_1$ and $\AAA_2$ are almost equal;
thus it is enough to find the limit of $\P\{\AAA_1\}$.
Indeed, let the indices $m$ and $n$ be such that $I_1=\{m, m+1, \ldots, n\}$. Then
\begin{eqnarray*}
	\P\{\AAA_2\} &\le& \P\{\AAA_1\} = \P\left\{ \bigcup_{i=m}^n \{ X_i^\eps \ge f((\ell+s)i)\eps \} \right\}
	\\
	&\le& \P\left\{X_m^\eps \ge f((\ell+s)m)\eps \right\} + \P\left\{X_{m+1}^\eps \ge f((\ell+s)(m+1))\eps \right\}
	\\
	&&  + \P\left\{ \bigcup_{i=m+2}^n \{ X_i^\eps \ge f((\ell+s)i)\eps \} \right\}
	\\
	&\le& 2\, \P\left\{X_m^\eps \ge f((\ell+s)m)\eps \right\}
	+ \P\left\{ \bigcup_{j=m+1}^{n-1} \{ X_{j+1}^\eps \ge f((\ell+s)(j+1))\eps \} \right\}
	\\
	&\le& 2\, \P\left\{X_m^\eps \ge f((\ell+s)m)\eps \right\} +   \P\{\AAA_2\},
\end{eqnarray*}
where in the penultimate inequality we use the stationarity of the sequence $X_i^\eps$ following
from the stationarity of the process $Y$.
For the remaining term we use the Pickands--Piterbarg asymptotics and obtain
\[
\P\left\{X_m^\eps \ge f((\ell+s)m)\eps \right\} \le \P\left\{X_m^\eps \ge f(\tau-\ell) \eps \right\}
\sim
c \ell (f(\tau-\ell)\eps)^{2/\alpha-1} \exp\{ -(f(\tau-\ell)\eps)^2/2\}.
\]

Using the relations \eqref{ls2}, \eqref{sigma0}, \eqref{ff'_asymp} we get that
$\ell = o(\sigma - \tau)$, $\ell f(\tau)f'(\tau) \eps^2\to 0$. Moreover, due to Lagrange's theorem and the relation \eqref{sigma0}, for some $\xi \in [\tau - \ell, \tau]$ one has 
\[
f(\tau - \ell)^2\eps^2 = (f(\tau)-\ell f'(\xi))^2 \eps^2 = f(\tau)^2\eps^2 - 2 \ell f(\tau)f'(\xi)\eps^2 + \ell^2 f'(\xi)^2 \eps^2 = f(\tau)^2\eps^2 + o(1)
\]
 Combining this with \eqref{tauassymp_1} and \eqref{tauassymp_0} we get that 
\begin{eqnarray*}
	\ell (f(\tau-\ell)\eps)^{2/\alpha-1} \exp\{ -(f(\tau-\ell)\eps)^2/2\}
	&\sim&
	\ell (f(\tau)\eps)^{2/\alpha-1} \exp\{ -(f(\tau)\eps)^2/2\}
	\\
	\sim
	\ell (f(\tau)\eps)^{1-\beta} e^{-r-\gamma^2/2}
	&\sim& \ell \gamma^{1-\beta} e^{-r-\gamma^2/2}.
\end{eqnarray*}
We know that
\[
\ell f'(\tau)f(\tau)\eps^2 \to 0,
\]
and from relation \eqref{tau_order}
\[
\ell f'(\tau)f(\tau)\eps^2 \asymp \ell\gamma^2\tau^{-\lambda}\asymp \ell \gamma^{1-\beta} e^{-\gamma^2/2}.
\]
That is why
\[
\ell \gamma^{1-\beta}e^{-r-\gamma^2/2} \to 0.
\]
We conclude that $\P\left\{X_m^\eps \ge(\ell+s)m\eps \right\}\to 0$, thus the difference between
$\P\{\AAA_1\}$ and $\P\{\AAA_2\}$ is indeed negligible.
\medskip

Below we repeatedly use the following technical lemma. Its proof is postponed to
Section \ref{ss:proof_prop}.
\begin{lem} \label{l:prop}
	For each $\alpha\neq 0$ and all $\theta(\eps), a(\eps), b(\eps), c(\eps)$ such that, as $\eps \to 0$, one has uniformly over $\omega\in [-1, 1]$
	$f'(\theta + \omega a)\sim f'(\theta)$, $f(\theta)\eps \sim \gamma$, $a = o(\theta)$, $f(\theta) c\eps^2\to 0$, $f(\theta)f'(\theta) a \eps^2\to 0$, it is true that 
	\begin{eqnarray*}
		\sum_{i: ai+b \ge\theta }^{\infty} (f(ai+b)\eps+c\eps)^{2/\alpha-1} \exp \{-(f(ai+b)\eps+c\eps)^2/2\}
		\\ \sim \frac{e^{\gamma^2/2}}{a} (f(\theta)\eps)^{2/\alpha+\beta-2}\exp\{-(f(\theta)\eps)^2/2\}.
	\end{eqnarray*}
\end{lem}

\bigskip

Let us evaluate $\P\{\AAA_3\}$. From stationarity of the sequence $V_i^{\eps}$ and Lemma \ref{1lem:pikpit} we obtain a uniform in $i$ bound
\begin{eqnarray*}
	\P\{\AAA_3\} &\le& \sum_{i\in I_3 } \P\{ V_i^\eps\ge f((\ell+s)i)\eps \}
	\\
	&\le& c\, s\, (1+o(1))\sum_{i: (\ell+s)(i+1)\geq \tau }    (f((\ell+s)i)\eps)^{2/\alpha-1}\exp \{- (f((\ell+s)i)\eps)^2/2\}.
\end{eqnarray*}

In order to find the asymptotic behavior of this sum, we apply Lemma \ref{l:prop} with parameters
$a=\ell+s, b=0, c=0, \theta = \tau-\ell-s$. Then, by using \eqref{ls1}, \eqref{ls2} and \eqref{tauassymp_1},
we have $a \sim \ell, \theta \sim \tau, (f(\theta) \eps)^2 = (f(\tau)\eps)^2+o(1)$.
Therefore, Lemma \ref{l:prop}, relation \eqref{tauassymp_0} and assumption \eqref{ls1} yield
\begin{eqnarray*}
	&&\frac{c s e^{\gamma^2/2} }{a}(f(\theta)\eps)^{2/\alpha+\beta-2}\exp\{-(f(\theta)\eps)^2/2\}
	\\&\sim&
	\frac{c s e^{\gamma^2/2}}{\ell}(f(\tau)\eps)^{2/\alpha+\beta-2}\exp\{-(f(\tau)\eps)^2/2\}
	= o(1).
\end{eqnarray*}

\medskip

Let us bound $\P\{\AAA_4\}$ as follows
\begin{eqnarray*}
	&& \P\{\AAA_4\} \le \sum_{j=0}^\infty  \P\{  \max_{t\in [ \sigma+ j, \sigma+j+1]} Y(t)>  f(\sigma + j)\eps \}
	\\
	&\le&  c\, (1+o(1)) \sum_{j=0}^\infty  (f(\sigma+ j)\eps)^{2/\alpha-1}\exp \{- (f(\sigma + j)\eps)^2/2\}.
\end{eqnarray*}

Lemma \ref{l:prop} applied with parameters $a=1, b=\sigma, c=0, \theta = \sigma$ and relation \eqref{sigmaassymp_0} provide the following asymptotics for the sum
\begin{eqnarray*}
	&&c e^{\gamma^2/2}(f(\sigma)\eps)^{2/\alpha+\beta-2}\exp\{-(f(\sigma)\eps)^2/2\} = o(1).
\end{eqnarray*}

\bigskip

The hardest part is to show that
\be
1- \P\{\AAA_1\}  \rightarrow \exp(-c\exp(-r)), \qquad \eps\to 0.
\ee

Our main tool here is the following classical inequality due to Slepian (see e.g.,
\cite[\S 14]{Lif_GSF}, \cite[Lecture 2]{Pit}).

\begin{lem}
	\label{Slepian_original}
	 Let $(U_1,...,U_n)$ and  $(V_1,...,V_n)$ be two centered Gaussian vectors such that
	$\E U_j^2=\E V_j^2$, $1\le j\le n$, and $\E (U_i U_j)\le \E (V_i V_j)$, $1\le i,j\le n$.
	Then for each $r\in \R$  one has
	\[
	\P\left\{ \max_{1\le j\le n} U_j \ge r \right\}
	\ge  \P\left\{ \max_{1\le j\le n} V_j \ge r \right\}.
	\]
\end{lem}

One may write this inequality in a slightly more general form (see \cite[Lecture 2]{Pit}): under assumptions of the
Slepian lemma, for all non negative $r_1,...,r_n$ one has
\[
\P\left\{ \exists j:\ U_j \ge r_j \right\}
\ge  \P\left\{ \exists j:\ V_j \ge r_j \right\}.
\]
This fact follows by application of the Slepian inequality to the vectors
$(\tfrac{U_1}{r_1},...,\tfrac{U_n}{r_n})$ and $(\tfrac{V_1}{r_1},...,\tfrac{V_n}{r_n})$
 and $r=1$.

The latter inequality easily extends to the Gaussian processes with continuous trajectories defined
on a metric space (by the way, the processes satisfying assumption \eqref{1assumption_0} belong
to this class). Namely, let $\{U(t),t\in T\}$ and $\{V(t),t\in T\}$
be two centered Gaussian processes with continuous trajectories defined on a common metric space $T$.
Let $\E U(t)^2=\E V(t)^2$, $t\in T$, and $\E (U(t_1) U(t_2))\le \E (V(t_1) V(t_2))$, $t_1,t_2\in T$.
Then for all compact sets $T_1,...,T_n$ in $T$ and all non negative $r_1,...,r_n$ it is true that
\be \label{1slepian}
\P\left\{ \bigcup_{j=1}^n  \left\{ \max_{t\in T_j} U(t) \ge r_j \right\} \right\}
\ge   \P\left\{  \bigcup_{j=1}^n  \left\{\max_{t\in T_j} V(t) \ge r_j\right\} \right\}.
\ee
\bigskip

Now we proceed to the proof of the remaining claim
\be
1- \P\{\AAA_1\}  \rightarrow \exp(-c\exp(-r)), \qquad \eps\to 0.
\ee
We provide the corresponding upper and lower bounds. In both cases we use the Slepian inequality in the form \eqref{1slepian}. Since $\AAA_1$ is defined by the process on long intervals $L_i, i\in I_1$, when reffering to long intervals in this part we mean $L_i, i\in I_1$.
\bigskip

\medskip

{\it Upper bound.}
Let us compare our process $Y$ with an auxiliary process $Z$ that is defined as follows.
First, let us consider a process $\widetilde{Y}(t), t \in \cup_{i\in I_1} L_i$
which consists of independent copies of $Y(t)$ on the intervals $L_i$.
Define
\[
\delta^2=\delta^2(\eps):= \sup_{t\ge s(\eps)} \left|\E( Y(t) Y(0) )\right|.
\]
Taking into account the correlation decay assumption \eqref{1assumption_infty} and assumption \eqref{ls0} concerning
the choice of $s$, we have, as $\eps\to 0$,
\begin{equation} \label{deltas}
	\delta^2= o((\ln s)^{-1}) = o(\gamma^{-2}).
\end{equation}
Let $\xi$ be an auxiliary standard normal random variable independent of the process $\widetilde{Y}$.
We define the centered Gaussian process $Z(t)$, $t \in \cup_{i\in I_1} L_i$, by the equality
\[
Z(t) := \sqrt{1-\delta^2} \widetilde{Y}(t) + \delta \xi.
\]
Then for all $t$ the variances are equal: $\E Y(t)^2= \E Z(t)^2=1$.
For covariances we have the following inequalities:
\begin{itemize}
	\item For $t_1$ and $t_2$ that belong to the same interval $L_i$ we have
	\begin{eqnarray*}
		\E(Z(t_1)Z(t_2))&=&
		\E\left(\left(\sqrt{1-\delta^2} \tY(t_1) + \delta \xi\right) \left(\sqrt{1-\delta^2} \tY(t_2)+\delta \xi\right)\right)
		\\&=&(1-\delta^2)\E(Y(t_1)Y(t_2))+\delta^2\ge \E(Y(t_1)Y(t_2)),
	\end{eqnarray*}
	where the last inequality follows from  $\E (Y(t_1)Y(t_2)) \le \sqrt{\E Y(t_1)^2 \E Y(t_2)^2}$ $=1$,
	\item for $t_1$ and $t_2$ that belong to different intervals $L_i$ and $L_j$,  by the definition of
	$\delta$ and by the intervals' construction we have
	\begin{eqnarray*}
		\E(Z(t_1)Z(t_2))&=&
		\E\left(\left(\sqrt{1-\delta^2} \tY(t_1) + \delta \xi\right) \left(\sqrt{1-\delta^2} \tY(t_2)+\delta \xi\right)\right)
		\\&=& \delta^2\ge \E(Y(t_1)Y(t_2)).
	\end{eqnarray*}
\end{itemize}

Let $\tX_i^\eps := \max_{t\in L_i} \tY(t)$. By applying Slepian inequality \eqref{1slepian}
to the processes $Y$ and $Z$, we obtain
\[
\P\{\AAA_1\}
= \P\left\{ \bigcup_{i\in I_1} \{ X_i^\eps \ge f((\ell+s)i)\eps \} \right\}
\ge \P\left\{ \bigcup_{i\in I_1} \{ \sqrt{1-\delta^2}\tX_i^\eps + \delta \xi \ge f((\ell+s)i)\eps \} \right\}.
\]
Let us pass to the complementary events; for every $h=h(\eps)>0$ the following elementary bound holds:
\begin{eqnarray} \nonumber
	1-\P\{\AAA_1\}
	&=& \P\left\{\bigcap_{i\in I_1} \{X_i^\eps\le f((\ell+s)i)\eps\}\right\}
	\\ \nonumber
	&\le& \P\left\{\bigcap_{i\in I_1} \{\sqrt{1-\delta^2} \tX_i^{\eps}+\delta \xi \le f((\ell+s)i)\eps\}\right\}
	\\ \nonumber
	&\le& \P\left\{\bigcap_{i\in I_1} \{\sqrt{1-\delta^2} \tX_i^{\eps}\le f((\ell+s)i)\eps + h\eps\}\right\}+\P\{\delta \xi \le -h\eps\}
	\\ \label{hepsdelta_prob}
	&=& \prod_{i\in I_1} \P\left\{ X_i^\eps\le \frac{f((\ell+s)i)\eps+h\eps}{\sqrt{1-\delta^2}} \right\}+\P\{\xi \le -h\eps/\delta\},
\end{eqnarray}
where the last equality holds because $\tX_i^\eps$  are independent copies of $X_i^\eps$.
We choose the level $h=h(\delta,\eps)$ so that, as $\eps\to 0$,
\begin{eqnarray} \label{h_large}
	&& h\eps/\delta \to \infty,
	\\  \label{h_small}
	&& h f(\tau)\eps^2 \sim h \gamma \eps \to 0,
\end{eqnarray}
which is possible under \eqref{deltas}.
\\Under \eqref{h_large} the last term of \eqref{hepsdelta_prob} is negligible.
Therefore, we aim to show that the product converges to $\exp(-c\exp(-r))$, as $\eps$ goes to 0.
Taking the logarithm and passing to the complementary events, we see that it is necessary to
prove the convergence
\[
\sum_{i\in I_1} \ln\left(1 - \P\left\{ X_i^\eps\ge \frac{f((\ell+s)i)\eps+h\eps}{\sqrt{1-\delta^2}} \right\}\right) \to - c\exp(-r).
\]
The probabilities in this sum tend to 0 uniformly over $i \in I_1$, since $X_i^{\eps}$ are identically distributed. Hence
\[
\sum_{i\in I_1} \ln\left(1 - \P\left\{ X_i^\eps\ge \frac{f((\ell+s)i)\eps+h\eps}{\sqrt{1-\delta^2}} \right\}\right) = -(1+o(1)) \sum_{i\in I_1} \P\left\{ X_i^\eps\ge \frac{f((\ell+s)i)\eps+h\eps}{\sqrt{1-\delta^2}} \right\}.
\]
Let us prove that 
\[
\sum_{i\in I_1} \P\left\{ X_i^\eps\ge \frac{f((\ell+s)i)\eps+h\eps}{\sqrt{1-\delta^2}} \right\} \to c\exp(-r).
\]
By Lemma \ref{1lem:pikpit} we know the exact asymptotics of each term of this sum, moreover due to stationarity of the sequence $X_i^{\eps}$ the equivalence is uniform over $i\in I_1$. Therefore the whole sum is equivalent to
\[
c \,\ell\, \sum_{i\in I_1} \left( \frac{f((\ell+s)i)\eps+h\eps}{\sqrt{1-\delta^2}} \right)^{2/\alpha-1}
\exp\left(-\left(\frac{f((\ell+s)i)\eps+h\eps}{\sqrt{1-\delta^2}} \right)^2/2\right).
\]
We represent this expression as the difference of two sums
\be \label{sum_1}
c \,\ell\, \sum_{i: (\ell+s)i+\ell\ge \tau} \left( \frac{f((\ell+s)i)\eps+h\eps}{\sqrt{1-\delta^2}} \right)^{2/\alpha-1}
\exp\left(-\left(\frac{f((\ell+s)i)\eps+h\eps}{\sqrt{1-\delta^2}} \right)^2/2\right)
\ee
and
\be \label{sum_2}
c \,\ell\, \sum_{i: (\ell+s) i \ge \sigma } \left( \frac{f((\ell+s)i)\eps+h\eps}{\sqrt{1-\delta^2}} \right)^{2/\alpha-1}
\exp\left(-\left(\frac{f((\ell+s)i)\eps+h\eps}{\sqrt{1-\delta^2}} \right)^2/2\right)
\ee
The asymptotics of the first sum follows from Lemma \ref{l:prop} applied with the variable change $\tilde{\eps}=\eps/\sqrt{1-\delta^2}\sim\eps$ and parameters
$a=\ell+s \sim \ell, b=0, c=h$ and $\theta=\tau-\ell \sim \tau$. Due to \eqref{h_small}, \eqref{tauassymp_1}, \eqref{ls2} and \eqref{deltas} it is true that
\begin{eqnarray*}
	(f(\theta)\tilde{\eps})^2&=& (f(\tau-\ell)\tilde{\eps})^2=(f(\tau)\tilde{\eps})^2+o(1)
	= (f(\tau)\eps)^2/(1-\delta^2) +o(1) 
	\\
	&=& (f(\tau)\eps)^2+O((f(\tau)\eps\delta)^2)+o(1) 
	= (f(\tau)\eps)^2+ o(1).
\end{eqnarray*}
Combining this with the condition \eqref{tauassymp_0}, Lemma \ref{l:prop}
gives us the asymptotics for \eqref{sum_1}
\begin{eqnarray*}
	&&c \,\ell\, \frac{e^{\gamma^2/2}}{a}(f(\theta) \tilde{\eps})^{2/\alpha+\beta-2}\exp\{-(f(\theta)\tilde{\eps} )^2/2\}
	\\ &\sim&  c \,\ell\, \frac{e^{\gamma^2/2}}{\ell} (f(\tau) \eps)^{2/\alpha+\beta-2}\exp\{-(f(\tau)\eps )^2/2\}
	\sim ce^{-r}.
\end{eqnarray*}
At the same time, Lemma \ref{l:prop} with the same $\tilde{\eps}$ and parameters $a=\ell+s\sim \ell, b=0, c=h$
and $\theta = \sigma$ provides us with asymtpotics of \eqref{sum_2}.
Since $(f(\theta)\tilde{\eps})^2 = (f(\sigma)\eps)^2+o(1)$ and \eqref{sigmaassymp_0}, we obtain
\begin{eqnarray*}
	&&c\,\ell\,\frac{e^{\gamma^2/2}}{a}(f(\theta)\tilde{\eps})^{2/\alpha+\beta-2}\exp\{-f(\theta\tilde{\eps})^2/2\}
	\\&\sim& c\,\ell\,\frac{e^{\gamma^2/2}}{\ell}(f(\sigma)\eps)^{2/\alpha+\beta-2}\exp\{-(f(\sigma)\eps)^2/2\}=o(1).
\end{eqnarray*}
Subtraction of the sums' asymptotics implies the required upper bound for $1-\P\{\AAA_1\}$.

\medskip

{\it Lower bound.}

In order to obtain an opposite bound for $1-\P\{\AAA_1\}$, we introduce and compare two more auxiliary processes
$Y_1$, $\tY_1$.
Let $\xi$ be an auxiliary standard normal random variable independent of the process $Y$.
Let $Y_1(t) := Y(t)+\delta\xi$, $t\in\cup_{i \in I_1} L_i$.  Furthermore, let us consider a sequence of independent standard Gaussian
random variables
$\xi_i$ independent of $\tY(t)$; and let
\[
\tY_1(t) := \tY(t)+\delta\xi_i, \qquad t\in L_i.
\]
Then for all $t$ we have the equality of variances: $\E Y_1(t)^2= \E\tY_1(t)= 1+\delta^2$.
For covariances we have the following inequalities:
\begin{itemize}
	\item For $t_1$ and $t_2$ that belong to the same interval $L_i$ we have
	\[
	\E(Y_1(t_1)Y_1(t_2)) = \E(\tY_1(t_1)\tY_1(t_2)),
	\]
	\item for $t_1$ and $t_2$ that belong to different intervals $L_i$ and $L_j$ we have
	\begin{eqnarray*}
		\E(Y_1(t_1)Y_1(t_2)) &=& \E((Y(t_1)+\delta\xi)(Y(t_2)+\delta\xi))=\E(Y(t_1)Y(t_2))+\delta^2
		\\
		&\ge& 0=	\E(\tY_1(t_1)\tY_1(t_2)).
	\end{eqnarray*}
\end{itemize}

We choose $h=h(\delta,\eps)$ as before, i.e. satisfying assumptions
\eqref{h_large} and \eqref{h_small}.

Slepian inequality \eqref{1slepian} yields
\begin{eqnarray*}
	&&  \P \left\{ \bigcup_{i\in I_1} \{ \tX_i^\eps+\delta \xi_i \ge f((\ell+s) i) \eps-h \eps \} \right\}
	= \P \left\{ \bigcup_{i\in I_1} \{ \max_{t\in L_i} \tY_1(t) \ge f((\ell+s) i) \eps-h \eps \} \right\}
	\\
	&\ge&  \P \left\{ \bigcup_{i\in I_1} \{ \max_{t\in L_i} Y_1(t) \ge f((\ell+s) i) \eps-h \eps \} \right\}
	=  \P \left\{ \bigcup_{i\in I_1} \{ X_i^\eps+\delta \xi \ge f((\ell+s) i) \eps-h \eps \} \right\}.
\end{eqnarray*}
By passing to the complementary events, we obtain
\begin{eqnarray*}
	&& \P \left\{ \bigcap_{i\in I_1} \{ X_i^\eps+\delta \xi \le f((\ell+s) i) \eps-h \eps \} \right\}
	\ge   \P \left\{ \bigcap_{i\in I_1} \{ \tX_i^\eps+\delta \xi_i \le f((\ell+s) i) \eps-h \eps \} \right\}
	\\
	&=& \prod_{i\in I_1}  \P  \left\{ \tX_i^\eps+\delta \xi_i \le f((\ell+s) i) \eps-h \eps \right\}
	=  \prod_{i\in I_1}  \P  \left\{ X_i^\eps+\delta \xi \le f((\ell+s) i) \eps-h \eps \right\}.
\end{eqnarray*}
Further, we apply an elementary bound
\begin{eqnarray*}
	1-\P\{\AAA_1\} &=& \P \left\{ \bigcap_{i\in I_1} \{ X_i^\eps \le f((\ell+s) i) \eps \} \right\}
	\\
	&\ge& \P \left\{ \bigcap_{i\in I_1} \{ X_i^\eps+\delta \xi \le f((\ell+s) i) \eps-h \eps \} \right\}  - \P \{ \delta \xi \le -h \eps \}
	\\
	&\ge&  \prod_{i\in I_1}  \P  \left\{ X_i^\eps+\delta \xi \le f((\ell+s) i) \eps-h \eps \right\} - \P \{ \delta \xi \le -h \eps \}.
\end{eqnarray*}
Under assumption \eqref{h_large}, we have $\P\{\delta \xi\le -h\eps\}\to 0$, as $\eps\to 0$.
It remains to prove that the product is greater than $\exp(-c\exp(-r))(1+o(1))$.
Taking the logarithm and passing to the complementary events, we see that it is necessary to
prove the bound
\[
\sum_{i\in I_1}\ln\left( 1 - \P \{X_i^\eps+\delta \xi \ge f((\ell+s) i) \eps-h \eps\}\right)
\ge -c\exp(-r) (1+o(1)).
\]
Since all the probabilities in this sum tend to zero uniformly over $I_1$, we have
\[
\sum_{i\in I_1}\ln\left( 1 - \P \{X_i^\eps+\delta \xi \ge f((\ell+s) i) \eps-h \eps\}\right)
=
- (1+o(1))\sum_{i\in I_1} \P \{X_i^\eps+\delta \xi \ge f((\ell+s) i) \eps-h \eps\}.
\]
Let us prove that
\[
\sum_{i\in I_1} \P \{X_i^\eps+\delta \xi \ge f((\ell+s) i) \eps-h \eps\}
\le c\exp(-r) (1+o(1)).
\]
We start with the estimate
\begin{eqnarray} \nonumber
	&& \sum_{i\in I_1} \P \{X_i^\eps+\delta \xi \ge f((\ell+s) i) \eps-h \eps\}
	\\
	\nonumber
	&\le& \sum_{i\in I_1} \big( \P \{X_i^\eps \ge f((\ell+s) i) \eps-2 h \eps\}
	+ \P\{\delta \xi>h\eps\}\big)
	\\ \label{N1}
	&\le&   \sum_{i: (\ell+s)i+\ell \geq \tau}  \P \{X_i^\eps \ge f((\ell+s) i) \eps-2 h \eps\}
	+ N_1 \P\{\delta \xi> h\eps\},
\end{eqnarray}
where $N_1$ denotes the number of elements in the set $I_1$; it has asymptotics
\[
N_1\sim \frac{\sigma-\tau}{\ell+s}=\frac{(R-r)B_\eps}{\ell+s}
\sim \frac{R}{f'(\tau)\eps\gamma \,\ell} \sim \frac{R}{f(\tau)f'(\tau)\eps^2\,\ell}.
\]
For the sum in \eqref{N1} Lemma \ref{1lem:pikpit} provides an equivalent expression
\be \label{sum_pickpit2}
c \,\ell\, \sum_{i\in I_1} \left( f((\ell+s) i) \eps-2h\eps \right)^{2/\alpha-1}
\exp\left(-\left(f((\ell+s) i) \eps -2 h\eps \right)^2/2\right).
\ee
Asymptotics for the last sum follows from Lemma \ref{l:prop} with parameters
$a=\ell+s$, $b=0$, $c=-2h$ $\theta = \tau-\ell$,
moreover, similarly to the upper bound, we have $a\sim\ell$, $\theta\sim\tau$, $(f(\theta)\eps)^2=(f(\tau)\eps)^2+o(1)$.
Therefore Lemma \ref{l:prop} combined with \eqref{tauassymp_0} provides us with a relation
\begin{eqnarray*}
	&&\frac{c\,\ell\, e^{\gamma^2/2}}{a} (f(\theta)\eps)^{2/\alpha+\beta-2} \exp\{-(f(\theta)\eps)^2/2\}
	\\&\sim& c\, e^{\gamma^2/2} (f(\tau)\eps)^{2/\alpha+\beta-2}\exp\{-(f(\tau)\eps)^2/2\}
	\sim ce^{-r}.
\end{eqnarray*}
We obtained the asymptotics of the first term in \eqref{N1}
\[
\sum_{i:(\ell+s)i+\ell \ge \tau}  \P \{X_i^\eps \ge f((\ell+s) i) \eps-2 h \eps\} = c \, e^{-r}(1+o(1)).
\]
It remains to estimate the second term in \eqref{N1}. To this end, we have to specify the choice of parameters
$\ell$ and $R$.

Since $\P(\delta\xi>h\eps)\to 0$, we may choose $\ell=\ell(\eps)$ that, although satisfying \eqref{ls2},
is still such that
\[
\frac{\P(\delta\xi>h\eps)}{\ell f(\tau)f'(\tau)\eps^2} \to 0.
\]
Then $R=R(\eps)$ can be choosen growing so slowly, that
\[
N_1 \, \P\{\delta \xi> h\eps\} \sim \frac{ R \, \P(\delta\xi>h\eps)}{\ell f(\tau)f'(\tau)\eps^2} \to 0.
\]
By summing up the estimates for the terms of \eqref{N1}, we arrive at the required lower estimate for $1-\P\{\AAA_1\}$.

\subsection{Proofs of technical lemmas} \label{ss:proof_prop}

\subsubsection{Proof of Lemma \ref{l:R}}
Let us remind the statement.

{\it	Consider any function $\tilde{R}(\eps)$ growing so slowly, that
	$\tilde{R} = o(\ln\gamma)$ and
	\[
	B_\eps \tilde{R} = \frac{\tilde{R}}{f'(\tau_0)\eps\gamma} + o\left(\frac{1}{f'(\tau_0)\eps\gamma}\right)
	\]
	Then for $A_\eps$ and $B_\eps$ described in Theorem $\ref{t:main}$ we have
	\begin{equation}
		f(A_\eps+ B_\eps \tilde{R})\eps = \gamma + \left(\frac{2}{\alpha}+\beta-2\right)\frac{\ln\gamma}{\gamma} + \frac{\tilde{R}}{\gamma} + o\left(\frac{1}{\gamma}\right).
	\end{equation}
	And, besides that, uniformly over $\omega\in [-1, 1]$ we have
	\[
	f'(\tau_0)\sim f'(\tau_0 + \omega (\sigma - \tau_0)), \;\;\;
	f(\tau_0) \sim f(\tau_0 + \omega (\sigma - \tau_0)).
	\]
}	

\begin{proof}

Define 
\[
\Delta = \frac{1}{f'(\tau_0)\eps}\left(\left(\frac{2}{\alpha}+\beta-2\right)\frac{\ln\gamma}{\gamma}
\\+\frac{\tilde{R}}{\gamma}\right).
\]
It's sufficient to show that
\[
f\left(\tau_0+\Delta + o\left(
\frac{1}{f'(\tau_0)\eps\gamma}\right)\right) = \frac{\gamma}{\eps} + \Delta f'(\tau_0)+o\left(\frac{1}{\gamma\eps}\right)=f(\tau_0)+\Delta f'(\tau_0)+o\left(\frac{1}{\gamma\eps}\right).
\]
Notice that
\[
\Delta \tau_0^{-\lambda}\sim \frac{\ln\gamma}{f'(\tau_0)\gamma\eps\tau_0^{\lambda}}
\\ \prec \frac{\ln\gamma}{f(\tau_0)\gamma\eps}\sim \frac{\ln\gamma}{\gamma^2}.
\]
In particular, $\Delta = o(\tau_0)$.
Moreover, uniformly over $\omega\in [-2, 2]$ we have
\[
f(\tau_0+\omega\Delta)\sim f(\tau_0) \;\;\;\mbox{and}\;\;\; f'(\tau_0+\omega\Delta)\sim f'(\tau_0),
\]
since for $\lambda < 1$
\begin{eqnarray*}
	\ln \left(\frac{f(\tau_0+\omega\Delta)}{f(\tau_0)}\right) = 
	\int_{\tau_0}^{\tau_0+\omega\Delta} (\ln f)'(x)\;dx \asymp \int_{\tau_0}^{\tau_0+\omega\Delta} x^{-\lambda}\;dx 
	\\ \asymp (\tau_0+\omega\Delta)^{1-\lambda}-\tau_0^{1-\lambda} \asymp \omega\Delta \tau_0^{-\lambda} = o(1),
\end{eqnarray*}
and for $\lambda=1$
\begin{eqnarray*}
	\ln \left(\frac{f(\tau_0+\omega\Delta)}{f(\tau_0)}\right) \asymp \int_{\tau_0}^{\tau_0+\omega\Delta} x^{-1}\;dx = \ln(1+\omega\Delta\tau_0^{-1})=o(1).
\end{eqnarray*}
The similar reasoning works for $f'$. This already shows the second part of the statement, since $\sigma - \tau_0$ can be represented as $\Delta$ for a right choice of $\tilde R$.
\\Uniformly over $\omega \in [0, 2]$ 
\[
f''(\tau_0+\omega\Delta)\asymp f'(\tau_0+\omega\Delta) (\tau_0+\omega\Delta)^{-\lambda}\sim f'(\tau_0)\tau_0^{-\lambda}\asymp f''(\tau_0). 
\]
By the mean value theorem for some $\omega\in [0, 2]$ we have
\[
f\left(\tau_0+\Delta+o\left(\frac{1}{f'(\tau_0)\eps\gamma}\right)\right)-f(\tau_0) = f'(\tau_0+\omega \Delta)\left(\Delta+o\left(\frac{1}{f'(\tau_0)\eps\gamma}\right)\right),
\]
since the additive in the argument is smaller $2\Delta$ as $\eps$ goes to 0.
Notice that
\[
f'(\tau_0+\omega\Delta)\cdot o\left(\frac{1}{f'(\tau_0)\eps\gamma}\right) = o\left(\frac{1}{\gamma\eps}\right).
\]
It remains to show that
\[
(f'(\tau_0+\omega\Delta)-f'(\tau_0))\Delta = o\left(\frac{1}{\gamma\eps}\right),
\]
or, equivalently,
\[
\frac{(f'(\tau_0+\omega\Delta)-f'(\tau_0))\ln\gamma}{f'(\tau_0)\eps\gamma} = o\left(\frac{1}{\gamma\eps}\right).
\]
By the mean value theorem for some $\tilde \omega\in[0, \omega]$
\[
f'(\tau_0+\omega\Delta)-f'(\tau_0) = \omega\Delta f''(\tau_0+\tilde\omega\Delta) \asymp \omega\Delta f''(\tau_0)\asymp \omega\Delta f'(\tau_0)\tau_0^{-\lambda}.
\]
We complete the proof by noticing that
\[
\Delta \tau_0^{-\lambda}\ln\gamma = o(1).
\]
\end{proof}
\subsubsection{Choice of $\ell(\eps)$ and $s(\eps)$}	

Let us remind that we are choosing $\ell$ and $s$ such that the relations \eqref{ls0}, \eqref{ls1}, \eqref{ls2} are true, i.e.
\begin{eqnarray*} 
	\ln s \succ \gamma^2, \mbox{i.e.} \ln s \succ -\ln(\eps) \; \mbox{for}\; \lambda=1 \mbox{ and }\ln s \succ \ln(-\ln\eps) \; \mbox{for}\; \lambda < 1,
	\\
	s/\ell\to 0,
	\\
	f(\tau) f'(\tau) \eps^2 \, \ell \to 0.
\end{eqnarray*}
We recall that due to \eqref{ff'_asymp} and \eqref{tau_order} it is true that
\[
f(\tau)f'(\tau)\eps^2 \asymp \gamma^2\tau^{-\lambda} \asymp \gamma^{1-\beta} e^{-\gamma^2/2} .
\]
Therefore we can pick $s = e^{\gamma^2/4}$ and any $\ell \succ c$ such that $\gamma^{1-\beta}e^{-\gamma^2/2}\ell \to 0$, as $\eps \to 0$. Note that $\ell = e^{\gamma^2/3}$ clearly fits these conditions, thus the set of choices for $\ell$ is not empty.

It is enough to show that for some $\delta>0$.
\[
e^{2\delta\gamma^2} \tau^{-\lambda} \to 0
\]
Then we can pick $s = e^{\delta\gamma^2/2}$ and any $\ell\succ s$ such that $\gamma^2 \tau^{-\lambda} \ell \to 0$. Note that $\ell = e^{\delta \gamma^2}$ fits these conditions so the set of choices for $\ell$ is not empty.

\subsubsection{Proof of lemma \ref{l:prop}}

Let us remind the statement.

\smallskip

	{\it For each $\alpha\neq 0$ and all $\theta(\eps), a(\eps), b(\eps), c(\eps)$ such that, as $\eps \to 0$, one has
	$f'(\theta +\omega a)\sim f'(\theta)$ and $f(\theta+\omega a)\sim f(\theta)$ uniformly over $\omega \in [-1, 1]$, $f(\theta)\eps \sim \gamma$, $a = o(\theta)$, $f(\theta) c\eps^2\to 0$, $f(\theta)f'(\theta) a \eps^2\to 0$, it is true that 
	\begin{eqnarray*}
		\sum_{i: ai+b \ge\theta }^{\infty} (f(ai+b)\eps+c\eps)^{2/\alpha-1} \exp \{-(f(ai+b)\eps+c\eps)^2/2\}
		\\ \sim \frac{e^{\gamma^2/2}}{a} (f(\theta)\eps)^{2/\alpha+\beta-2}\exp\{-(f(\theta)\eps)^2/2\}.
	\end{eqnarray*} }

\begin{proof}

The monotone decay of the function $x\mapsto x^{2/\alpha-1}\exp \{- x^2/2\}$ at large $x$ combined with $(f(\theta-a)+c)\eps \to \infty$ allows to bound the sum between two integrals:
\[
\frac{1}{a}\int_{\theta\pm a}^{\infty}  (f(x)\eps+c\eps)^{2/\alpha-1}\exp \{- (f(x)\eps+c\eps)^2/2\} dx. 
\]
By changing the variables $y=(f(x)+c)\eps$ we get $x=f^{-1}(y/\eps-c)$. Then $dx = \frac{1}{\eps}(f^{-1})'(y/\eps-c)dy$. The integrals transform into
\[
\frac{1}{a\eps}\int_{(f(\theta\pm a)+c)\eps}^{\infty}  y^{2/\alpha-1}\exp \{- y^2/2\} (f^{-1})'(y/\eps-c)dy.
\]
To deal with these integrals we estimate $(f^{-1})'(y/\eps-c)$ for large $y$ and for small $y$ separately via regularity conditions \eqref{reg1} and \eqref{reg2} and the relation.
\[
(f(\theta\pm a)+c)\eps \sim \gamma.
\]
Since $f(\theta\pm a)\eps\sim\gamma$, as $\eps$ goes to 0
\[
[f(\theta\pm a)\eps, (1+\kappa/2)f(\theta\pm a)\eps] \subset [(1-\kappa)\gamma, (1+\kappa)\gamma].
\]
Therefore, regularity condition \eqref{reg1} gives us
\begin{eqnarray*}
	&&\frac{1}{a\eps}\int_{(f(\theta\pm a)+c)\eps}^{(1+\kappa/2)(f(\theta\pm a)+c)\eps} y^{2/\alpha-1}\exp\{-y^2/2\}(f^{-1})'(y/\eps-c) dy
	\\ &\sim& \frac{(f^{-1})'(1/\eps)}{a\eps}\int_{(f(\theta\pm a)+c)\eps}^{(1+\kappa/2)(f(\theta\pm a)+c)\eps} y^{2/\alpha-1}\exp\{-y^2/2\}(y-c\eps)^{\beta} dy
	\\ &\sim& \frac{(f^{-1})'(1/\eps)}{a\eps}\int_{(f(\theta\pm a)+c)\eps}^{(1+\kappa/2)(f(\theta\pm a)+c)\eps} y^{2/\alpha+\beta-1}\exp\{-y^2/2\} dy
	\\ &\sim& \frac{(f^{-1})'(1/\eps)}{a\eps}((f(\theta\pm a)+c)\eps)^{2/\alpha+\beta-2}\exp\{-((f(\theta\pm a)+c)\eps)^2/2\}
	\\ &\sim& \frac{(f^{-1})'(1/\eps)}{a\eps}(f(\theta)\eps)^{2/\alpha+\beta-2}\exp\{-(f(\theta)\eps)^2/2\}.
\end{eqnarray*}
The last equivalence follows from the Lagrange's theorem and our conditions, since for some $\omega\in [-1, 1]$ one has 
\begin{eqnarray*}
(f(\theta\pm a) + c)^2\eps^2 &=& (f(\theta) + f'(\theta + \omega a) a + c)^2\eps^2
\\&=& f(\theta)^2\eps^2 + O(f(\theta) f'(\theta+\omega a) a \eps^2 + f(\theta) c \eps^2) 
= f(\theta)^2\eps^2 + o(1).
\end{eqnarray*}
Besides that, the second regularity condition $\eqref{reg2}$ provides us with
\[
(f^{-1})'\left(\frac{y-c\eps}{\eps}\right) < (y-c\eps)^{\tilde\beta}(f^{-1})'\left(\frac{1}{\eps}\right),
\]
followed by
\begin{eqnarray*}
	\frac{1}{a\eps}\int_{(1+\kappa/2) (f(\theta\pm a)+c)\eps}^{\infty} y^{2/\alpha-1}\exp\{-y^2/2\}(f^{-1})'(y/\eps-c)dy
	\\ \prec \frac{(f^{-1})'(1/\eps)}{a\eps} ((1+\kappa/2) (f(\theta\pm a)+c)\eps)^{2/\alpha+\tilde \beta-2}\exp\{-((1+\kappa/2) (f(\theta\pm a)+c)\eps)^2/2\}
	\\=o\left(\frac{(f^{-1})'(1/\eps)}{a\eps}(f(\theta)\eps)^{2/\alpha+\beta-2}\exp\{-(f(\theta)\eps)^2/2\}\right).
\end{eqnarray*}
Finally, we get that the integral asymptotics is
\[
\frac{(f^{-1})'(1/\eps)}{a\eps}(f(\theta)\eps)^{2/\alpha+\beta-2}\exp\{-(f(\theta)\eps)^2/2\} \sim  \frac{e^{\gamma^2/2}}{a} (f(\theta)\eps)^{2/\alpha+\beta-2}\exp\{-(f(\theta)\eps)^2/2\}.
\]

\end{proof}

\medskip

\bigskip


\end{document}